%% file: ccp.tex
%

\magnification 1200
\input amstex
\documentstyle{amsppt}
\rightheadtext{Cardinal Characteristics and $\Pi$}

\define\p{\hbox{$\Pi$}}
\define\z{\hbox{$\Bbb Z$}}
\define\s{\hbox{$\Sigma$}}
\define\add{\bold{add}}
\define\cov{\bold{cov}}
\define\unif{\bold{unif}}
\define\pp{{\frak p}}
\define\dd{{\frak d}}
\define\bb{{\frak b}}
\define\cc{{\frak c}}
\define\se{{\frak {se}}}
\define\ee{{\frak e}}
\define\tl{torsionless}

\topmatter
\title
Cardinal Characteristics and the \\
Product of Countably Many\\
Infinite Cyclic Groups
\endtitle
\author
Andreas Blass
\endauthor
\address
Mathematics Dept., University of Michigan, Ann Arbor, MI
48109, U.S.A.
\endaddress
\email
ablass\@math.lsa.umich.edu
\endemail
\thanks Partially supported by NSF grant DMS-9204276.
\endthanks
\abstract
We study, from a set-theoretic point of view, those subgroups of
the infinite direct product $\z^{\aleph_0}$ for which all
homomorphisms to \z\ annihilate all but finitely many of the
standard unit vectors.   Specifically, we relate the smallest
possible size of such a subgroup to several of the standard
cardinal characteristics of the continuum.  We also study some
related properties and cardinals, both group-theoretic and
set-theoretic. One of the set-theoretic properties and the
associated cardinal are combinatorially natural, independently
of any connection with algebra.
\endabstract
\endtopmatter
\document

\head
Introduction
\endhead
Let $\p=\z^{\aleph_0}$  be the direct product of a countable
infinity of copies of the infinite cyclic group \z.
Specker [19] proved that \p\ as well as many of its subgroups
$G$
have the following property, in which $e_n$ is the element of
\p\ whose $n$th component is 1 and whose other components
are all zero (the $n$th standard unit vector).  If $h$ is a
homomorphism from \p\ (or $G$) to \z, then $h(e_n)=0$ for all
but finitely many $n$. The subgroups $G$ for which Specker
established this property all have, like \p\ itself, the cardinality
of the continuum, ${\frak c}=2^{\aleph_0}$, and the question
naturally arises whether any smaller subgroups $G$ of \p\ also
have this property.  Eda [7] showed that this question is
undecidable on the basis of the usual (Zermelo-Fraenkel) axioms
of set theory (ZFC).  Specifically, he proved that the answer is
negative in models of Martin's Axiom  but positive in models
obtained by adjoining many random reals.  In fact, his proofs
give somewhat more precise information about the minimum
cardinality $\kappa$ of a subgroup of \p\ satisfying Specker's
theorem.  In terms of the cardinal characteristics of the
continuum introduced in [6] and described in Section 1 below,
Eda's proofs establish that
${\frak p}\leq\kappa\leq{\frak d}$.

One purpose of this paper is to improve these estimates; we
shall show that the additivity of Lebesgue measure
${\bold{add}}(L)$
 is another lower bound for $\kappa$ and that
the bounding number $\frak b$  is an
upper bound.
(Definitions of these cardinals are recalled in Section 1.)
The new upper bound subsumes Eda's, since
${\frak b}\leq{\frak d}$; the new lower bound is
incomparable with Eda's, since neither of $\frak p$ and
${\bold{add}}(L)$ is provably larger than the other.

These proofs suggest some additional group-theoretic questions
concerning homomorphisms into free groups.  We say that an
abelian group $G$ ``binds'' a subgroup $H$ if every
homomorphism of $G$ into a free group maps $H$ to a group of
finite rank.  We consider questions about the smallest
cardinalities in which various sorts of binding (a group binding
an infinite rank subgroup, a group binding itself, etc.) can occur
non-trivially.   A rather surprising result is that one of these
cardinalities is $\geq{\bold{add}}(L)$ but $\leq$ the additivity
of Baire category, ${\bold{add}}(B)$.
These questions are connected to the Specker phenomenon
described above and also to questions considered by Eklof and
Shelah [9] about groups that satisfy
$G\cong G\oplus F$ with $F$ free of finite rank.

Finally, we consider some purely set-theoretic questions arising
out of these problems.  These questions seem quite natural from
a combinatorial point of view, but do not seem to have been
previously considered.

This paper is organized as follows.  Section 1
presents the definitions of the cardinal characteristics of the
continuum that we will need later and some relevant known
theorems about them.  Section 2 is devoted to the Specker
phenomenon and the proof that it occurs in cardinality
$\frak b$, i.e., the upper bound mentioned above.  The lower
bound is established in Section 4 as a consequence of a stronger
result about binding.  Section 3 begins with a general discussion
of binding,describes its connection with the Specker
phenomenon and with the cardinal $\pp$,
and ends with the connection with the Eklof-Shelah work
mentioned above.  Finally, Section 4 contains the purely
set-theoretic notions of ``predicting'' and ``evading,'' their
connection with cardinal characteristics and applications to the
group-theoretic topics of the previous sections.

\head
Acknowledgement
\endhead

I thank John Irwin for many helpful conversations about
abelian groups and for his tireless efforts to keep my attention
focused on \p.

\head
Terminology and Notation
\endhead

By ``group'' we mean abelian group, and we write the group
operation as addition.
\z\ is the group of integers.
For any set $I$, $\z^I$ is the group of functions $I\to\z$ with
addition of corresponding values as the operation; it is the
direct product of $|I|$ copies of \z, where $|I|$ means the
cardinality of $I$.
$\z^{(I)}$ is the subgroup of $\z^I$ consisting of functions whose
values are 0 at all but finitely many elements of $I$; it is the
direct sum of $|I|$ copies of \z.
When $I$ is the set $\omega$ of natural numbers, we write
\p\ and \s\ for $\z^\omega$ and $\z^{(\omega)}$, respectively.
If $x$ and $y$ are in $\z^I$ and at least one of them is in
$\z^{(I)}$, then the sum $\sum_{i\in I}x(i)y(i)$ is finite, and we
denote it by the inner product notation $\langle x,y\rangle$.

A group $G$ is {\sl \tl\/} if for every non-zero $x\in G$ there is
a homomorphism $h:G\to\z$ with $h(x)\neq0$.  This is
equivalent to requiring that $G$ be embeddable in $\z^I$ for
some set $I$. (For one direction of the equivalence, use the
projection homomorphisms $\z^I\to\z$ that evaluate functions
in $\z^I$ at a specific element of $I$.  For the other direction,
take $I$ to be the set of all homomorphisms $G\to\z$ and
embed $G$ in $\z^I$ by $g\mapsto(i\mapsto i(g))$.)  In any
group $G$, the elements that are mapped to 0 by all
homomorphisms $G\to\z$ form a subgroup $N$, and the
quotient $G/N$ is the largest \tl\ quotient of $G$.

We use the standard (among set-theorists) notation $\omega$
for the set of natural numbers.  When discussing functions on
$\omega$ or subsets of $\omega$, we often use an asterisk~$^*$
to indicate that finitely many exceptions are allowed.  For
example, if $A$ and $B$ are subsets of $\omega$, then
$A\subseteq^*B$ means that $A\diagdown B$ is finite.
Similarly, if $f$ and $g$ are functions from $\omega$ into an
ordered set, then $f\leq^*g$ means that $f(n)\leq g(n)$ for all
but finitely many $n$.

\input ccp1
\input ccp2

\input ccp3

\input ccp4

\head
5. Questions
\endhead

The results presented in this paper raise a multitude of
questions, for we have introduced numerous cardinals and
many possible connections between them remain undecided.
Here are some rather strong conjectures, whose proof would
greatly simplify our picture of these cardinals.
\roster
\item $\ee-=\se$. (Notice that this would imply equality of
cardinals \therosteritem2 through \therosteritem5 in
Corollary~8.)
\item $\ee\leq\bb$.
\item $\se=\add(B)$.
\endroster
A more specific problem, whose solution might throw
considerable light on the general situation, is to compute the
value of $\se$ in the model obtained by an $\aleph_2$-stage,
countable support iteration of Mathias forcing over a model of
the generalized continuum hypothesis.  This model has $\bb$
and both $\unif$'s equal to $\aleph_2$ while both $\cov$'s are
$\aleph_1$.  Therefore $\ee=\aleph_1$ by Theorem~13, but our
results do not determine $\se$.

\input ccpref

\enddocument

%% file: ccp1.tex
%

\head
1. Cardinal Characteristics of the Continuum
\endhead

For information about cardinal characteristics of the continuum,
including the facts stated without proof below, we refer to
Vaughan [20] and the earlier works cited there.
We shall have occasion to refer to nine of the standard cardinal
characteristics of the continuum, three related to Lebesgue
measure, three related to Baire category, and three of a more
combinatorial nature.

The characteristics related to Lebesgue measure are
\roster
\item $\add(L)$, the additivity of measure, i.e., the smallest
number of measure zero sets whose union is not of measure
zero,
\item $\cov(L)$, the covering number for measure, i.e., the
smallest number of measure zero sets that can cover the real
line, and
\item $\unif(L)$, the uniformity number for measure, i.e., the
smallest cardinality of a set not of measure zero.
\endroster
To avoid possible confusion, we emphasize that ``not of measure
zero'' is not synonymous with ``of non-zero measure,'' because
the former includes the possibility of not being measurable.

The analogous characteristics for Baire category are
\roster\item[4] $\add(B)$, the smallest number of first category
sets whose union is of second category,
\item $\cov(B)$, the smallest number of first category sets that
can cover the real line, and
\item $\unif(B)$, the smallest cardinality of a set of second
category.
\endroster

These definitions referred to measure and category in the real
line, but the cardinals defined here would be the same if we
used either of the following two spaces instead.  $2^\omega$
is the set of infinite sequences of zeros and ones, with the
product topology obtained from the discrete topology on
$\{0,1\}$ and the product measure obtained from the measure
on $\{0,1\}$ that gives each of the two elements measure
$1/2$.  Similarly $\omega^\omega$ is the set of infinite
sequences of non-negative integers, with the product topology
obtained from the discrete topology on $\omega$ and with the
product measure obtained from the measure on $\omega$ that
gives each point $n$ measure $2^{-(n+1)}$.

To define the three remaining characteristics that we will need,
we recall the notations $\subseteq^*$ and $\leq^*$ introduced
above, and we also introduce the following terminology.  A
family $\Cal F$ of sets is said to have the {\sl strong finite
intersection property\/} if the intersection of every finite
subfamily of $\Cal F$ is infinite.  In terms of these concepts we
define
\roster
\item[7] $\pp$ is the smallest cardinality of a family $\Cal F$ of
subsets of $\omega$ such that $\Cal F$ has the strong finite
intersection property but there is no infinite
$A\subseteq\omega$ satisfying $A\subseteq^*F$ for all
$F\in\Cal F$,
\item $\dd$ is the smallest cardinality of a family $\Cal D$ of
functions from $\omega$ to $\omega$ such that every function
from $\omega$ to $\omega$ is $\leq^*$ some function in
$\Cal D$, and
\item $\bb$ is the smallest cardinality of a family $\Cal B$ of
functions from $\omega$ to $\omega$ such that no function
from $\omega$ to $\omega$ is $\geq^*$ all the functions in
$\Cal B$.  (In other words, for every $g:\omega\to\omega$
there is $f\in\Cal B$ such that $g(n)<f(n)$ for infinitely
many~$n$.)
\endroster

Clearly, $\pp$ would be unchanged if $\omega$ in its definition
were replaced by any countably infinite set.  Similarly, $\dd$
and $\bb$ would be unchanged if their definitions referred to
functions $A\to\omega$ instead of $\omega\to\omega$, where
$A$ is countably infinite.  Also, $\dd$ would be unchanged if
the relation of ``almost everywhere majorization,'' $\leq^*$,
were replaced with ``everywhere majorization,'' for we can
adjoin to $\Cal D$, without increasing its cardinality, all those
functions that differ only finitely from functions in $\Cal D$.

For orientation, we point out that all nine of these cardinals lie
between $\aleph_1$ and the cardinality $\cc$ of the continuum,
inclusive.  In particular, if the continuum hypothesis holds, then
all are equal to $\aleph_1=\cc$.  It follows from results of
Martin and Solovay [14] that under Martin's axiom all these
cardinals equal $\cc$.

We shall need some information about the relationships
between these nine cardinals.  Most of this information is
summarized in the following diagram, in which an arrow from
one cardinal to another indicates that the former is provably
less than or equal to the latter; that is $\rightarrow$ means
$\leq$.  (Except for the part involving $\pp$, this is part of what
is called Cicho\'n's diagram; see [3, 10, 20].  The use of arrows
instead of inequality signs is due to the custom of working not
with the cardinals $\kappa$ themselves but with the
hypotheses $\kappa=\cc$; inequalities between the cardinals
obviously yield implications between these hypotheses, and in
practice the converse also holds, i.e., the known proofs of
implications between these hypotheses also establish the
inequalities.)

$$\matrix
\cov(L)&\longrightarrow&\unif(B)&&&&\\
\Big\uparrow&&\Big\uparrow&&&&\\
\smash{\big\vert}&&\bb&\longrightarrow&\dd&&\\
\Big\vert&&\Big\uparrow&&\Big\uparrow&&\\
\add(L)&\longrightarrow&\add(B)&\longrightarrow&\cov(B)&
     \longrightarrow&\unif(L)\\
&&\Big\uparrow&&&&\\
&&\pp&&&&
\endmatrix$$

A few of the inequalities are obvious from the definitions,
namely $\add\leq\cov$ and $\add\leq\unif$ for both measure
and category, and $\bb\leq\dd$. The inequalities
$\cov(B)\leq\dd$ and $\bb\leq\unif(B)$ follow easily from the
observation that, for each $g\in{}^\omega\omega$, the set of all
$f\leq^*g$ is of first category in ${}^\omega\omega$.

Among the non-trivial inequalities, $\cov(L)\leq\unif(B)$ and
$\cov(B)\leq\unif(L)$ are due to Rothberger [18],
$\add(L)\leq\add(B)$ is due to Bartoszy\'nski [1] and
independently to Raisonnier and Stern [17], $\pp\leq\add(B)$ is
due to Martin and Solovay [14],
and $\add(B)\leq\bb$ is due to Miller [16].  In fact, Miller
proved a stronger result, which we shall need later and
therefore display here for emphasis:
$$
\add(B)=\min\{\cov(B),\bb\}.\tag10
$$

The diagram above is complete in the sense that every provable
(in ZFC) inequality between our nine cardinals is given by a
chain of arrows in the diagram.  For inequalities not involving
$\pp$, the proofs are summarized in [3].  The consistency of
$\pp<\add(L)$ is proved in [13], and the consistency of
$\pp>\cov(L)$ is stated in [11], attributed to Miller
(unpublished).

To close this section, we retract our earlier statement that
$\pp$, $\dd$, and $\bb$ are more combinatorial than the
measure and category characteristics.  Not that the former
aren't combinatorial, but it turns out that the latter have
combinatorial descriptions as well.  We shall need these
descriptions, due to Bartoszy\'nski [1,2], for $\add(L)$ and
$\cov(B)$; for related results, see also [15,16, 17].

For the first of these descriptions, we use Bartoszy\'nski's
notion of a {\sl slalom\/}, namely a function $s$ with domain
$\omega$, such that $s(n)$ is a set of cardinality $(n+1)^2$ for
each $n$.  We say that a function $f:\omega\to\omega$
{\sl goes through\/} the slalom $s$ if $f(n)\in s(n)$ for all but
finitely many $n$. With this terminology, we have that
$\add(L)$ is the smallest possible cardinality of a family $\Cal
F$ of functions  $\omega\to\omega$ such that, for each slalom
$s$, some $f\in\Cal F$ does not go through $s$.

We say that two functions $f$ and $g$ on $\omega$ are
infinitely often equal if $f(n)=g(n)$ for infinitely many $n$.
Then $\cov(B)$ is the smallest possible cardinality of a family
$\Cal F$ of functions  $\omega\to\omega$ such that no single
function $g:\omega\to\omega$ is infinitely often equal to each
of the functions in $\Cal F$.

%% file: ccp2.tex
%

\head
2. The Specker Phenomenon
\endhead

We say that a subgroup $G$ of $\p=\z^{\omega}$
{\sl exhibits the Specker phenomenon\/} if it contains a
sequence $(a_n)_{n\in\omega}$ of linearly independent
elements such that every homomorphism $G\to\z$ maps $a_n$
to 0 for all but finitely many $n$.  (The concept, but not the
name, is from [7].) When we need to be more specific, we shall
say that $G$ exhibits the Specker phenomenon {\sl witnessed
by\/} $(a_n)_{n\in\omega}$.
Thus, Specker's theorem [19, Satz III] asserts that \p\ exhibits
the Specker phenomenon witnessed by the sequence of
standard unit vectors $(e_n)$.

\definition{Definition}
The smallest cardinal of any subgroup of \p\ exhibiting the
Specker phenomenon is denoted by $\se$.
\enddefinition

The symbol $\se$ stands for Specker and Eda.  Eda studied this
cardinal in [7], proving that $\se=\cc$ follows from Martin's
axiom but is false in a model obtained by adding many random
reals.  In fact, his argument for the random real case shows that
$\se\leq\dd$.  (It is well-known that $\dd<\cc$ in the random
real model.)  His proof that Martin's axiom implies $\se=\cc$
also establishes somewhat more, though indirectly.  The
application of Martin's axiom in the proof uses a partially
ordered set which does not merely (as Martin's axiom requires)
satisfy the countable antichain condition but in fact is
$\sigma$-centered.  Bell [4] showed that Martin's axiom for
$\sigma$-centered orderings is equivalent to $\pp=\cc$.  In
conjunction with Eda's argument, this proves that $\pp=\cc$
implies $\se=\cc$, and a trivial variation of the argument shows
that $\pp\leq\se$.  (We shall give a simpler proof of this, not
using Bell's theorem but following the lines of Eda's proof, in
Section 3.) Thus, we can summarize Eda's results (in the light of
Bell's theorem) as $\pp\leq\se\leq\dd$.

In this section, we shall improve the upper bound on $\se$
from $\dd$ to $\bb$.  Later, in Section 4, we shall give a new
lower bound, namely $\add(L)$.

Before proceeding, it seems appropriate to make some remarks
about the definition of the Specker phenomenon.  On the one
hand, the definition might seem too restrictive --- why are only
subgroups of \p\ considered?  On the other hand, one might be
even more restrictive and bring the definition closer to
Specker's result by requiring that the witnesses $a_n$ be the
standard unit vectors $e_n$.  We wish to point out that, once
some trivialities are removed, such variations in the definition
do not affect $\se$.  As an example of a triviality to be
removed, we note that, if $G$ were not required to be a
subgroup of \p, then every divisible group of infinite rank
would exhibit the Specker phenomenon simply because it has
no non-zero homomorphisms to \z.

Quite generally, if $G$ is not \tl, then its homomorphisms to \z\
are ``the same'' as those of its largest \tl\ quotient $G/N$, so,
when dealing with homomorphisms to \z\ (e.g., when discussing
the Specker phenomenon), it is reasonable to divide out the
irrelevant $N$ and work with the \tl\ $G/N$.  Thus, when
relaxing the restriction, in the definition of Specker
phenomenon, that $G$ be a subgroup of \p, we should still
require that it be torsionless.

\proclaim{Theorem 1}
Either of the following two modifications of the definition of
$\se$ does not change its value.\newline
(a) Instead of requiring $G$ to be a subgroup of\/ \p, require
only that $G$ be \tl.\newline
(b) Require the witness sequence to be the sequence $(e_n)$ of
all standard unit vectors.
\endproclaim

\demo{Proof}
(a)\enspace The assumption that $G$ is \tl\ is equivalent to
requiring that $G$ be embeddable in $\z^I$ for some index set
$I$.  Suppose now that $G\subseteq\z^I$ exhibits the Specker
phenomenon, witnessed by $(a_n)$.  For each non-zero linear
combination $c$ of (any finitely many) $a_n$'s, select one
coordinate $i\in I$ such that the $i$th component of $c$ is not
zero.  Let $J$ be the set of $i$'s so selected; note that $J$ is a
countable subset of $I$.  Let $p:\z^I\to\z^J$ be the canonical
projection map.  By our choice of $J$, the elements $p(a_n)$ are,
like the $a_n$, linearly independent.  For any homomorphism
$h:p(G)\to\z$, we know that $hp:G\to\z$ annihilates almost all
$a_n$, so $h$ annihilates almost all $p(a_n)$.  Thus, $p(G)$
exhibits the Specker phenomenon witnessed by $(p(a_n))$.
Since $J$ is countable, $p(G)$ is isomorphic to a subgroup of \p,
so we have an instance of the Specker phenomenon as originally
defined, involving a group $p(G)$ no larger than the given $G$.
Thus, $\se$ defined with the relaxed condition on $G$ in (a) is
no smaller than $\se$ as originally defined.  That it is no larger
is trivial, so (a) is proved.

(b)\enspace Consider any $G\subseteq\p$ exhibiting the
Specker phenomenon witnessed by $(a_n)$.  The first part of
the proof of Specker's [19, Satz I] can be used to show that,
given any countably many elements of \p, like the $a_n$, we
have an endomorphism of \p\ mapping all these elements into
the subgroup \s\ of elements having only finitely many
non-zero components.  (See also Chase [5, Corollary 3.3].)  Thus,
we may assume without loss of generality, replacing $G$ by its
image under a suitable endomorphism of \p, that our witnesses
$a_n$ all lie in \s.  Let us write $S_n$ for the (finite) set of
coordinates where $a_n$ has non-zero components.  Each
coordinate $i$ is in only finitely many $S_n$'s, because the
projection of $G\subseteq\p$ to the $i$th coordinate is a
homomorphism to \z, to which the definition of the Specker
phenomenon is applicable.  This observation allows us to find an
infinite subsequence of $(a_n)$ for which the supports $S_n$
are pairwise disjoint.  Indeed, we can put $a_0$ into our
subsequence, put out those finitely many other $a_n$ whose
$S_n$ meets $S_0$, put in the first $a_k$ different from all
these, put out the finitely many other $a_n$ whose $S_n$ meets
$S_k$, etc.  Of course, the resulting subsequence, like any
infinite subsequence of $(a_n)$, still witnesses the Specker
phenomenon for $G$, so we may assume, without loss of
generality, that $(a_n)$ is such that the supports $S_n$ are
pairwise disjoint.  By applying suitable automorphisms of \p,
composed from automorphisms of all the finite rank free groups
$\z^{S_n}$ separately, we can arrange that each $a_n$ is an
integer multiple of a standard unit vector $e_n$.  Replacing $G$
by its smallest  supergroup pure in \p\ (which does not
increase the cardinality), we can arrange that the Specker
phenomenon is witnessed by a sequence of (some of the)
standard unit vectors, say the $e_n$ for $n$ in a certain infinite
subset $T$ of $\omega$.  We can arrange further that $G$ is a
subgroup of $\z^T$, just by replacing it with its image under the
projection $\p=\z^\omega \to \z^T$.  Finally, using a bijection
between $T$ and $\omega$, we obtain a subgroup of \p, of no
greater cardinality than the original $G$, for which the Specker
property is witnessed by the sequence of all the standard unit
vectors $e_n$.  This proves (b).
\qed\enddemo

The preceding theorem shows that the definition of $\se$ is
robust under fairly wide variations in the sort of group and the
sort of witnessing sequence considered.  Now we turn to our
first estimate of $\se$; the proof uses ideas from Eda [8].

\proclaim{Theorem 2}
\p\ has a subgroup of cardinality $\bb$ that exhibits the
Specker phenomenon (witnessed by the sequence $(e_n)$ of
standard unit vectors).  Thus $\se\leq\bb$.
\endproclaim

\demo{Proof}
Let $\Cal B$ be a family of $\bb$ functions $\omega\to\omega$
as in the definition \thetag9 of $\bb$, i.e., for any
$g:\omega\to\omega$ there is $f\in\Cal B$ such that $g(n)<f(n)$
for infinitely many $n$.  We can arrange that the functions in
$\Cal B$ are monotone non-decreasing and nowhere zero; just
replace each $f\in\Cal B$ with the larger function
$f'(n)=1+\max\{f(k)\mid k\leq n\}$. Henceforth, we assume
that this replacement has been made, so $0<f(m)\leq f(n)$
whenever $m\leq n$ and $f\in\Cal B$.

\proclaim{Lemma}
For every infinite $A\subseteq\omega$ and every function
$g:A\to\omega$, there is $f\in\Cal B$ such that $g(n)<f(n)$ for
infinitely many $n\in A$.
\endproclaim

\demo{Proof of lemma}
Extend the given $g$ to a function on all of $\omega$ by
defining its value at any $k\notin A$ to be the (already
defined) value of $g$ at the next larger element of $A$. By our
choice of $\Cal B$, it contains an $f$ that is greater than (the
extended) $g$ at infinitely many $k$.  If infinitely many of
these $k$'s are in $A$, we are done.  Otherwise, for each of
these $k$'s that is not in $A$, let $n$ be the next larger element
of $A$, and observe that, by definition of $g(k)$ and
monotonicity of $f$, we have $g(n)=g(k)<f(k)\leq f(n)$. So $f$ is
greater than $g$ at infinitely many $n\in A$.
\qed\enddemo

For each $f\in\Cal B$, define a function $x:\omega\to\omega$
by the recursion $x(0)=1$ and
$$
x(n+1)=2\cdot f(n+1)\cdot x(n)\cdot\sum_{i=0}^n x(i).\tag{11}
$$
(The motivation for this strange definition will become clear in
the course of the proof.)
We note for future reference that $x(n)$ is never zero, that
$x(n)\to\infty$ as $n\to\infty$, and that, because of the factor
$x(n)$ in the definition, $x(n)$ divides $x(m)$ whenever $n\leq
m$.

As $\Cal B$ has cardinality $\bb$, we have obtained $\bb$
functions $x$ in this manner.  Let $G$ be the smallest pure
subgroup of \p\ containing all these $x$'s and containing the
(countably many) elements of \s\ (the functions in \p\ with
finite support).  Then $G$ has cardinality $\bb$, and we
complete the proof by showing that $G$ exhibits the Specker
phenomenon witnessed by $(e_n)$.

We proceed by contradiction, so suppose $h:G\to\omega$ is a
homomorphism and the set  $A=\{n\mid h(e_n)\neq0\}$ is
infinite.  Define a function $g$ by $g(n)=\max_{k\leq n}
|h(e_n)|$, and apply the lemma to obtain an $f\in\Cal B$ such
that
$M=\{n\in A\mid f(n)>g(n)\}$ is infinite.  Let $x$ be the
function defined, using this $f$, by \thetag{11}. Thus $x\in G$.
Since $x(n)\to\infty$ as $n\to\infty$ and since $M$ is infinite,
we can fix a non-zero $n\in M$ with $x(n)>2|h(x)|$.  We shall
prove that $h(e_n)=0$; this will be the desired contradiction,
because $n\in M\subseteq A$.

We split $x\in G$ as the sum
$$
x=\sum_{i<n}x(i)e_i+y=\sum_{i<n}x(i)e_i+x(n)z,
$$
where $y$ is the element of \p\ that agrees with $x$ from
coordinate $n$ on but is zero at the earlier coordinates.  Since
all components $x(m)$ for $m\geq n$ are divisible by $x(n)$,
we can write $y$ as $x(n)z$ for some $z\in\p$. Furthermore, as
$G$ contains $x$ and $\sum_{i<n}x(i)e_i$
(the latter because it is in \s), it also contains their difference
$y$ and, being pure in \p, also contains $z$. Thus, we have
$$
h(x)=\sum_{i<n}x(i)h(e_i)+x(n)h(z)\equiv\sum_{i<n}x(i)h(e_i)
\pmod{x(n)}.\tag{12}
$$
Recall that, by our choice of $n$, the left side $h(x)$ of this
congruence has absolute value smaller than half the modulus.
So does the right side, because
$$
\left|\sum_{i<n}x(i)h(e_i)\right|\leq\sum_{i<n}x(i)g(n)
<f(n)\sum_{i<n}x(i)\leq\frac{1}{2}x(n),\tag{13}
$$
where we have used the definition of $g$, the fact that $n\in
M$, and the definition \thetag{11} of $x$.
But two numbers congruent modulo $x(n)$ and each smaller in
absolute value than half of that modulus must be equal. So we
have
$$
h(x)=\sum_{i<n}x(i)h(e_i).\tag{14}
$$

The preceding calculation can be re-done with $n+1$ in place of
$n$.  Some care is needed because $n+1$, unlike $n$, need not
be in $M$. This affects only the calculation \thetag{13}, which
becomes
$$
\left|\sum_{i<n+1}x(i)h(e_i)\right|\leq\sum_{i\leq n}x(i)g(n)
<f(n)\sum_{i\leq n}x(i)\leq f(n+1)\sum_{i\leq n}x(i)
\leq\frac{1}{2}x(n+1).
$$
Then we obtain, by the same argument as for \thetag{14},
$$
h(x)=\sum_{i<n+1}x(i)h(e_i).
$$
Subtracting \thetag{14} from this, we obtain $x(n)h(e_n)=0$.
Since $x$ is nowhere zero, we conclude that $h(e_n)=0$, the
desired contradiction to $n\in A$.
\qed\enddemo

%% file: ccp3.tex
%

\head
3. Binding Subgroups
\endhead

\definition{Definition}
A group $G$ {\sl binds\/} a subgroup $H$ if every
homomorphism from $G$ to a free group maps $H$ into a group
of finite rank.
\enddefinition

This definition has some trivial cases, e.g., if $H$ has finite rank
or if $G$ admits no non-zero homomorphisms to free groups.
As in the discussion preceding Theorem~1, elements of $G$
whose image under every homomorphism to \z\ is 0 are
irrelevant distractions in connection with binding, so it is
reasonable to divide $G$ by the subgroup of these elements, i.e.,
to replace $G$ by its largest \tl\ quotient.
Thus, when we discuss binding, we shall always assume that
$G$ is torsionless and therefore embeddable in $\z^I$ for some
set $I$.  We shall also assume, to avoid triviality, that $H$ has
infinite rank.

Notice that, if $G$ binds $H$, then every (\tl) supergroup of $G$
binds every (infinite rank) subgroup of $H$ as well.

\definition{Definition}
If a \tl\ group $G$ binds some subgroup of infinite rank, then
we say $G$ is {\sl binding\/}.  If $G$ has infinite rank and binds
itself, then we call it {\sl self-binding}.
\enddefinition

Clearly, every self-binding group is binding.

\proclaim{Proposition 3}
A \tl\ group $G$ binds a subgroup $H$ if and only if every
homomorphism $G\to\s$ maps $H$ to a group of finite rank.
$G$ binds itself if and only if it does not admit a homomorphism
onto \s, if and only if it does not have \s\ as a direct summand.
\endproclaim

\demo{Proof}
For the non-trivial half of the first statement, suppose
$f:G\to\z^{(I)}$ maps $H$ to a group of infinite rank.  We must
achieve the same situation with a countable set in place of $I$.
As $f(H)$ has infinite rank, we can choose a countable infinity
of linearly independent elements in it and then choose, for each
non-zero linear combination $x$ of these elements, one element
$i\in I$ with $x(i)\neq0$.  Let $J$ consist of the countably many
$i$ so chosen, and compose $f$ with the canonical projection
$\z^{(I)}\to\z^{(J)}$.  The resulting homomorphism
$G\to\z^{(J)}\cong\s$
is as desired.

For the non-trivial half of the second statement, suppose $G$
does not bind itself, and choose, by what we just proved,
$f:G\to\s$ with $f(G)$ of infinite rank.  But any subgroup of \s\
of infinite rank is isomorphic to \s, so we can compose $f$ with
an isomorphism $f(G)\to\s$ to obtain a surjection as claimed.
The last part of the proposition, about direct summands, follows
because any homomorphism onto a free group, such as \s, splits.
\qed\enddemo

We shall need to consider groups $G$ with $\s\subseteq
G\subseteq\p$.  We call such a group {\sl\s-binding\/} if it
binds \s, and we call it {\sl weakly \s-binding\/} if every
homomorphism $G\to\s$ that extends to an endomorphism of
\p\ maps \s\ to a group of finite rank.  It is clear, by the first
part of Proposition 3, that \s-binding implies weakly
\s-binding.  Although \s-binding is a special case of the general
notion of binding, the following theorem shows that it is not
very special.

\proclaim{Theorem 4}
If a \tl\ group $G$ binds a subgroup of infinite rank, then $G$
has a homomorphic image $G'\subseteq\p$ such that $G'+\s$
binds \s.
\endproclaim

\demo{Proof}
Let $G\subseteq\z^I$ bind a subgroup $H$ of infinite rank.  As
$G$ binds all
subgroups of $H$, we may assume that $H$ is countable.
Choosing for each
non-zero element $x\in H$ one $i\in I$ such that $x(i)\neq0$,
letting
$J\subseteq I$ be the set of these countably many chosen $i$'s,
and
replacing $G$ and $H$ by their projections in $\z^J$, we arrange
that, up to
isomorphism, $G\subseteq\p$, $G$ still binds $H$, and $H$ still
has infinite
rank.

The next step in the proof resembles an argument in Specker's
proof of [19, Satz I]
and even more closely resembles Theorem 3.2 of [5].

\proclaim{Lemma}
If $H$ is any subgroup of \p\ of countably infinite rank, then
there exists an endomorphism $f$ of \p\ and there exist
positive integers $(d_n)_{n\in\omega}$ such that $f(H)$
includes the subgroup $\bigoplus_n d_n \z$ of \s.
\endproclaim

\demo{Proof of lemma}
We define the endomorphism $f$ as the composition of an
infinite sequence of endomorphisms $f=\dots\circ f_2\circ
f_1\circ f_0$; the composition will be well-defined because
$f_k$ will not change the $n$th component of $x$ for any
$x\in\p$ and any $n<k$.  Thus the $n$th component of $f(x)$ is
the $n$th component of
$f_n\circ\dots\circ f_2\circ f_1\circ f_0(x)$.

To define $f_0$, choose a non-zero element $a\in H$ such that
the greatest common divisor $d$ of all its components $a(n)$ is
as small as possible.  Notice that $d$ is the largest integer by
which $a$ is divisible in \p. This provides a description of $d$
that is invariant under automorphisms of \p.  Also, $d$ is an
integral linear combination of finitely many $a(n)$'s, say
$\sum_{n<r}c_na(n)$.  Let $g$ be the automorphism of \p\ that
leaves the first $r$ components of its argument unchanged but
transforms the rest according to
$$
g(x)(m)=x(m)-\frac{a(m)}{d} \sum_{n<r}c_n x(n) \quad
\text{for }m\geq r.
$$
Recall that $d$ divides all $a(m)$, so this makes sense, and
notice that $g(a)(m)=0$ for all $m\geq r$.  Now, by applying a
suitable automorphism of $\z^r$ to the first $r$ coordinates
(where $g(a)$ can be non-zero) and leaving all the other
coordinates unchanged, we can send $g(a)$ to $d\cdot e_0$.
(The components of $g(a)$ have greatest common divisor $d$,
thanks to the automorphism-invariance of $d$ noted above.)

Let $f_0$ be the composite of $g$ and the automorphism just
described.  So $f_0(a)=d\cdot e_0$.  Also, let $d_0=d$.

Consider any element $b\in f_0(H)$.  As $f_0$ is an
automorphism of \p, the greatest common divisor of the
components of $b$ is the same as for
$f_0^{-1}(b)\in H$, hence is at least $d$.  In particular, $|b(0)|$
is at least $d$, unless it is zero.  It follows that $b(0)$ is
divisible by $d$, for otherwise we could write $b(0)=qd+r$ with
$0<r<d$, and then $b'=b-qde_0\in f_0(H)$ would have 0th
component $r$, a contradiction because the conclusion of the
previous sentence applies to $b'$ as well as to $b$.

Thus, every element of $f_0(H)$ is expressible as $d_0e_0+x$,
where $x\in f_0(H)$ has $x(0)=0$. Therefore, $f_0(H)=\langle
d_0e_0\rangle\oplus H_1$, where $H_1$ is a subgroup of
$\p_1=\z^{\omega\diagdown\{0\}}$.

Now proceed with $H_1$ just as we did with $H$ in the
preceding paragraphs, obtaining an automorphism $f_1$ of
$\p_1$ (which we extend to \p\ by letting it act trivially on the
0 coordinate) sending $H_1$ to
$\langle d_1e_1\rangle\oplus H_2$, where $H_2$ is a subgroup
of $\p_2=\z^{\omega\diagdown\{0,1\}}$ and where $d_1$ is the
smallest g.c.d. of all the components of any element of $H_1$.

Continuing this process inductively, we obtain $f_n$'s whose
composite sends $H$ to a subgroup containing all the $d_ne_n$.
(Notice that, although each $f_n$ is an automorphism of \p, we
can only claim that the composite $f$ is an endomorphism.)
This completes the proof of the lemma.
\qed\enddemo

Returning to the proof of the theorem, we apply the lemma to
the $H$ that was bound by our (modified) $G\subseteq\p$, and
we set $G'=f(G)$, where $f$ is given by the lemma.  Then $G'$ is
a homomorphic image of our modified $G$, hence also of our
original $G$.  It binds $f(H)$, so its supergroup $G'+\s$ binds the
subgroup $S=\bigoplus_n d_n \z$ of $f(H)$.
But $S$ is a subgroup of \s\ such that the quotient $\s/S$ is a
torsion group, because $S$ contains a multiple of each of the
generators $e_n$ of \s.  It follows immediately that any
homomorphism that is defined on \s\ and maps $S$ into a
group of finite rank must also map \s\ into a group of (the
same) finite rank.  Thus, binding $S$ implies binding \s, and the
proof is therefore complete.
\qed\enddemo

\proclaim{Corollary 5}
Each of the following cardinals is less than or equal to the next:
\roster
\item $\aleph_1$
\item the smallest cardinality of a weakly \s-binding group
\item the smallest cardinality of a binding group
\item the smallest cardinality of a self-binding group
\item the cardinality $\cc$ of the continuum.
\endroster
\endproclaim

\demo{Proof}
\thetag1$\leq$\thetag2:\enspace
We must show that, if $\s\subseteq G\subseteq\p$ and $G$ is
countable, then there is an endomorphism of \p\ mapping $G$
into \s\ and mapping \s\ onto a group of infinite rank.  For this,
it more than suffices to get an automorphism of \p\ mapping
$G$ into \s. Furthermore, it suffices to do this under the
additional assumption that $G$ is a pure subgroup of \p,
because $G$ can be enlarged to a pure subgroup without
increasing its cardinality.  But now what we need is given by
[5], Corollary 3.3.

\thetag2$\leq$\thetag3:\enspace
By Theorem 4, the cardinal in \thetag3 is also the smallest
cardinality of a \s-binding group.  Since \s-binding implies
weakly \s-binding, the desired inequality follows.

\thetag3$\leq$\thetag4:\enspace
This is immediate, as every self-binding group is binding.

\thetag4$\leq$\thetag5:\enspace
\p\ has cardinality $\cc$ and binds itself. (If it had \s\ as a
direct summand, then its dual, which is isomorphic to \s\ by
[19, Satz III], [12, Cor. 94.6],
would have a summand isomorphic to the dual of \s, i.e., to \p.
This is absurd as \s\ is countable and \p\ is  not.  The same
conclusion also follows instantly from the more detailed
information in [12, pp.159--160] about homomorphisms from
products to sums.)
\qed\enddemo

The next two results relate binding to the Specker phenomenon
discussed in Section 2.

\proclaim{Theorem 6}
A group that exhibits the Specker phenomenon is binding.
\endproclaim

\demo{Proof}
Suppose $G$ exhibits the Specker phenomenon witnessed by
$(a_n)$, and let $H$ be the subgroup of $G$ generated by these
witnesses $a_n$.  Suppose also, toward a contradiction, that $G$
does not bind $H$.  So let $f:G\to\s$ map $H$ to a group of
infinite rank.  We inductively construct an element $x\in\p$
and an increasing sequence of natural numbers $k_m$ as
follows.

At stage $m$ of the construction, we have already defined $k_i$
for $i<m$, and we have already defined some finite initial
segment of $x$.  As the group $f(H)$  generated by the
$f(a_n)$'s has infinite rank, we can choose $k_m$ so that
$f(a_{k_m})$ has a non-zero component in some position, say
the $i$th, such that $x(i)$ is not yet defined.  Then we can easily
extend $x$ (or, rather, the finite part of $x$ already defined) so
that (a)~$x(j)$ becomes defined for all $j$ such that
$f(a_{k_m})(j)\neq0$ and so that (b)~the inner product
$\langle f(a_{k_m}),x\rangle$ is not zero.  (The inner product
will not depend on future steps in the definition of $x$, because
of (a).  Making it non-zero is easy, by appropriately choosing
the value of
$x(i)$.)

Notice that all the $k_m$ are distinct, because the $m$th stage
is the first one where $x$ is defined at all coordinates where
$f(a_{k_m})$ has a non-zero component, so we can recover $m$
from $k_m$.

Now the homomorphism $z\mapsto\langle z,x\rangle$ maps \s\
and hence $f(G)$ into \z\ and maps none of the $f(a_{k_m})$'s
to zero.  So its composite with $f$ violates the choice of the
$a_n$ as witnessing the Specker phenomenon.
\qed\enddemo

This theorem immediately implies that $\se$ is at least as big as
the cardinal labeled \thetag3 in the Corollary 5.  In fact, we can
do a bit better, replacing \thetag3 with \thetag4.

\proclaim{Theorem  7}
The smallest cardinality of a self-binding group is $\leq\se$.
\endproclaim

Thus, $\se$ could be inserted into Corollary 5 as item
\thetag{4.5}.  Theorems 6 and 7 cannot be combined to assert
that a group exhibiting the Specker phenomenon binds itself; a
counterexample is given by $\s\oplus\p$ regarded as a
subgroup of \p\ in the obvious way (sequences in which only
finitely many odd-numbered components are non-zero).

\demo{Proof of Theorem 7}
By Theorem 1, let $G_1$ be a group of cardinality $\se$ such
that $\s\subseteq G_1\subseteq\p$ and $G_1$ exhibits the
Specker phenomenon witnessed by the sequence of standard
unit vectors $e_n$.  By Theorem 6, $G_1$ is binding, and by
Corollary 5 there is a weakly \s-binding group $G_2$ of
cardinality at most $\se$.  Let $G_3=G_1+G_2$.  Then
$\s\subseteq G_3\subseteq\p$, $G_3$ exhibits the Specker
phenomenon witnessed by the $e_n$'s, $G_3$ is weakly
\s-binding, and $|G_3|=\se$.

The next step of the proof is based on an idea from [8].  Fix two
disjoint increasing sequences of prime numbers, say $(p_n)$
and $(q_n)$.  For each $n$, obtain from the Chinese remainder
theorem an integer $d_n$ that is divisible by all $p_k$ for
$k\leq n$ but is congruent to $-1$ modulo all $q_k$ for $k\leq
n$.  Extend $G_3$ to a pure subgroup $G$ of \p\ that is closed
under (component-by-component) multiplication by the
sequence $d=(d_n)$.  This can be done with $|G|=\se$,
because we are simply closing $G_3$ under countably many
partial functions (multiplication by $d$ and division by each
positive integer).  The following lemma is due to Eda [8].

\proclaim{Lemma}
A homomorphism $f:G\to\z$ is completely determined if the
values $f(e_n)$ are known.
\endproclaim

\demo{Proof}
It suffices to show that, if $f(e_n)=0$ for all $n$ then $f(x)=0$
for all $x\in G$.  Fix any such $f$ and $x$.  Also, temporarily fix
a positive integer $n$. Decompose the componentwise product
$dx$ as the sum of $\sum_{i<n}d_ix(i)e_i$ and the rest, which,
by choice of $d_n$, is divisible by $p_n$ in \p, hence also in $G$
by purity.  The first part of this decomposition is annihilated by
$f$, as all $f(e_i)=0$. So $f(dx)$ is divisible by $p_n$.
Now un-fix $n$.  $f(dx)$ is divisible by arbitrarily large primes
$p_n$, hence is zero.  We can repeat the same argument with
the sequence $d+1=(d_n+1)$ in place of $d$ and with $q_n$ in
place of $p_n$ (since $q_k$ divides $d_n+1$ when $k\leq n$).
We obtain $f((d+1)x)=0$.  But then, as $x=(d+1)x-dx$, it follows
that $f(x)=0$, as desired.
\qed\enddemo

Combining the lemma with the Specker phenomenon witnessed
by the $e_n$'s, we find that every homomorphism $f:G\to\z$
has the form $f(x)=\sum_{i<n}f(e_i)x(i)$ for some $n$.
Therefore $f$ extends to a homomorphism $\p\to\z$.  It
follows, by applying this observation to each component, that
every homomorphism $G\to\s$ extends to a homomorphism
$\p\to\p$.
We use this to complete the proof of the theorem by showing
that $G$ binds itself.

Let $f:G\to\s$.
As $f$ extends to a homomorphism $\p\to\p$ and as $G$ is
weakly \s-binding, $f(\s)$ must have finite rank.  Therefore, for
all but finitely many $n\in\omega$, the $n$th component
$f_n:G\to\z$ of $f$ is zero on \s. But, by the lemma, this implies
that these (all but finitely many) $f_n$ are zero on $G$.  So
$f(G)$ has finite rank, as required.
\qed\enddemo

The following corollary extends Corollary 5 to incorporate
Theorems 2 and 7.

\proclaim{Corollary  8}
Each of the following cardinals is less than or equal to the next:
\roster
\item $\aleph_1$
\item the smallest cardinality of a weakly \s-binding group
\item the smallest cardinality of a binding group
\item the smallest cardinality of a self-binding group
\item $\se$
\item $\bb$
\item $\cc$.\qed
\endroster
\endproclaim

We are now in a position to fulfill our promise, from Section 2,
to give another proof, not using Bell's theorem, of the inequality
$\pp\leq\se$.  In fact, we obtain a stronger result, showing that
$\pp\leq$ all the cardinals except $\aleph_1$ mentioned in
Corollary~8.

\proclaim{Theorem  9}
A group of cardinality $<\pp$ cannot be weakly \s-binding.
\endproclaim

\demo{Proof}
Suppose $\s\subseteq G\subseteq\p$ and $|G|<\pp$.  Recall
from Section 1 that if $C$ is a countable set, if $\Cal F$ is a
family of fewer than $\pp$ subsets of $C$, and if every finite
subfamily of $\Cal F$ has infinite intersection, then there is an
infinite subset $A$ of $C$ such that $A\subseteq^*F$ for all
$F\in\Cal F$.  We apply this with $C=\s\diagdown\{0\}$ and
with the following $|G|$ sets $Z_x$ as the family $\Cal F$.  For
each $x\in G$, let
$$
Z_x=\{c\in C\mid \langle c,x\rangle=0\}.
$$

To check that the family $\Cal F=\{Z_x\mid x\in G\}$ has the
strong finite intersection property, let any finitely many of its
elements, say $Z_{x_1}$, $Z_{x_2}$, \dots, $Z_{x_n}$, be given;
we must show that their intersection is infinite.
For each $k\in\omega$, let $y_k$ be the $n$-component integer
vector $(x_1(k), x_2(k), \dots, x_n(k))$.  These infinitely many
vectors $y_k$, lying in the finite-dimensional vector space
$\Bbb Q^n$, must be linearly dependent; fix a non-zero
$r$-tuple $(c(0), c(1),\dots,c(r))$ of rational numbers such that
$\sum_{k=0}^r c(k)x_i(k)=0$ for all $i=1,2,\dots, n$.  Clearing
denominators, we can arrange that the $c(k)$ are integers.
Extending the definition of $c(k)$ to $k>r$ by making these
$c(k)=0$, we obtain an element $c\in\s\diagdown\{0\}$ with
$\langle c,x_i\rangle=0$ for all $i=1,2,\dots, n$.  This shows that
the intersection of the $Z_{x_i}$ is non-empty; in fact, it is
infinite because we can multiply $c$ by any non-zero integer.

Having checked the strong finite intersection property, we use
the fact that $|\Cal F|\leq|G|<\pp$ to obtain an infinite set
$A=\{a_0,a_1,\dots,a_n,\dots\}\subseteq\s\diagdown\{0\}$
such that $A\subseteq^*Z_x$ for all $x\in G$.
Now we define a homomorphism $f:\p\to\p$ by $f(x)(n) =
\langle a_n,x\rangle$.

For each $x\in G$, we have, for all but finitely many
$n\in\omega$, that $a_n\in Z_x$, which means that $f(x)(n)=0$.
Thus, $f(x)\in\s$.  We have shown that (the restriction of) $f$
maps $G$ into \s\ and extends to an endomorphism (namely
the unrestricted $f$) of \p.  We shall show that $f(\s)$ has
infinite rank, thereby completing the proof that $G$ is not
weakly \s-binding.

Suppose, toward a contradiction, that $f(\s)$ had finite rank.
Then there would be an $m\in\omega$ such that, for all $n\geq
m$ and all $k\in\omega$, $f(e_k)(n)=0$.  This equation means
$0=\langle a_n,e_k\rangle=a_n(k)$.  This would make all the
$a_n$ for $n>m$ equal to zero, contradicting the fact that the
$a_n$ are all in $C$ which doesn't contain zero.
\qed\enddemo

This theorem allows us to insert $\pp$ as item \thetag{1.5} in
Corollaries 5 and 8.

We close this section by connecting the notion of binding to a
question studied by Eklof and Shelah in [9].  One of their results
(not the main one) is that, under the assumption of Martin's
axiom, if $G$ is a group of cardinality $<\cc$ and if $G\cong
G\oplus\z^m$ for some positive integer $m$, then $G\cong
G\oplus\s$ (and therefore $G\cong G\oplus\z^m$ for every
positive integer $m$).  In fact, their application of Martin's
axiom involved a
$\sigma$-centered partially ordered set, so, by Bell's theorem,
the hypotheses that Martin's axiom holds and $|G|<\cc$ can be
weakened to the hypothesis that $|G|<\pp$.  Here is an alternate
approach, not using Bell's theorem, to (a stronger form of) this
result.

Assume $G\cong G\oplus\z^m$ for some positive integer $m$.
Instead of assuming $|G|<\pp$, we assume only that $|G|$ is
smaller than the smallest cardinality of any self-binding group.
(This hypothesis would follow if $|G|<\pp$, by virtue of
Corollary~8 and Theorem~9.)  Let $N$ be, as before, the
subgroup of $G$ consisting of elements mapped to 0 by all
homomorphisms $G\to\z$.  The analogous subgroup of
$G\oplus\z^m$ is clearly $N\oplus0$, as $\z^m$ is \tl.  So the
assumed isomorphism from $G$ onto $G\oplus\z^m$ must send
$N$ to $N\oplus0$ and must therefore induce an isomorphism
from the largest \tl\ quotient $G/N$ onto $(G/N)\oplus\z^m$.
By comparing the ranks of these groups, we conclude that
$G/N$ has infinite rank.  (This is the only use we need to make
of the assumption that $G\cong G\oplus\z^m$.)  As
$|G/N|\leq|G|<$ the smallest cardinal of any self-binding group,
Proposition 3 tells us that $G/N$ admits a surjection to \s, and
therefore so does $G$, and therefore $G\cong A\oplus\s$ for
some $A$.  But since $\s\cong\s\oplus\s$, it follows that
$G\cong G\oplus\s$, as desired.

%% file: ccp4.tex
%

\head
4. Predicting and Evading
\endhead

\define\rt{\restriction}

This section is devoted to combinatorial concepts, prediction and
evasion, motivated by some of the group-theoretic concepts of
the preceding section, particularly the weakest of these, weak
\s-binding.

\definition{Definition}
For any set $S$, an {\sl $S$-valued predictor\/} is a pair
$\pi=(D_\pi,(\pi_n)_{n\in D_\pi})$
where $D_\pi$ is an infinite subset of $\omega$ and $\pi_n$ is,
for each $n\in D_\pi$, a function $S^n\to S$.  We say that the
predictor $\pi$ {\sl predicts\/} the function
$x:\omega\to S$ if, for all but finitely many $n\in D_\pi$, we
have $x(n)=\pi_n(x\rt n)$.  Otherwise, we say that $x$ {\sl
evades\/} $\pi$.
\enddefinition

Here $x\rt n$ means the $n$-tuple $(x(0),x(1),\dots,x(n-1))$.
Intuitively, we regard predictors $\pi$ as follows.  The values of
an unknown function $x:\omega\to S$ are being revealed, one
at a time, in order, and we are to guess $x(n)$ just before it is to
be revealed, i.e., just after we have seen $x\rt n$.
The predictor $\pi$ provides a strategy for making these
guesses when $n\in D_\pi$, namely, if we have seen the
$n$-tuple $t$ so far, then we are to guess $\pi_n(t)$.  The
functions predicted by $\pi$ are just those for which all but
finitely many of the guesses provided by this strategy are
correct.

\definition{Definition}
A \z-valued predictor $\pi$ is called {\sl linear\/} if, for each
$n\in D_\pi$, the function $\pi_n:\z^n\to\z$ is a linear function
with rational coefficients.
A predictor is {\sl forgetful\/} if, whenever $m<n$ are
consecutive elements of $D_\pi$, the function $\pi_n$ depends
only on the last $n-m-1$ components of its argument $n$-tuple.
\enddefinition

In terms of the intuitive picture of predictors, forgetfulness
means that, when guessing $x(n)$, the strategy considers only
the values of $x$ revealed since its last guess, namely $x(m+1)$
through $x(n-1)$.

\definition{Definition}
$\ee$, the {\sl evasion number\/}, is the smallest possible
cardinality of a family $\Cal E$ of functions $\omega\to\omega$
such that every
$\omega$-valued predictor is evaded by some $x\in\Cal E$.
$\ee_l$ (resp. $\ee_f$, resp. $\ee_{fl}$)  is the smallest possible
cardinality of a family $\Cal E$ of functions $\omega\to\z$ such
that every linear (resp. forgetful, resp. forgetful linear)
$\z$-valued predictor is evaded by some $x\in\Cal E$.
\enddefinition

Note that $\ee$ would be unchanged if its definition were
phrased in terms of functions $\omega\to D$ and $D$-valued
predictors, where $D$ is any countably infinite set.  In
particular, $D$ could be \z, and so it is clear that $\ee$ is
greater than or equal to both $\ee_l$ and $\ee_f$, which are in
turn greater than or equal to $\ee_{fl}$.

The general notion of predictor and the associated cardinal
$\ee$ seem quite natural, but it is the more specialized notions
involving linearity and the associated cardinals that connect
directly with the group-theoretic concepts of the preceding
section, as the following theorem shows.

\proclaim{Theorem 10}
The following three cardinals are equal.
\roster
\item The smallest cardinality of a weakly \s-binding group
\item $\ee_l$
\item $\ee_{fl}$
\endroster
\endproclaim

\demo{Proof}
\therosteritem1$\geq$\therosteritem2:\enspace
Let $G$ be a weakly \s-binding group.  In particular,
$G\subseteq\p$, so $G$ is a family of functions $\omega\to\z$.
We shall show that every linear predictor $\pi$ is evaded by
some element of $G$, so $\ee_l\leq|G|$, as required.

Suppose, toward a contradiction, that $\pi$ is linear and
predicts every element of $G$.  Let $D_\pi=\{d_0,d_1,\dots\}$.
Thus, for each $x\in G$, all but finitely many $k\in\omega$
have $x(d_k)=\pi_{d_k}(x\rt d_k)$.  As $\pi_{d_k}$ is linear
with rational coefficients, we can clear denominators to rewrite
this equation as a linear relation with integer coefficients
$$
\sum_{i=0}^{d_k} c_{ki}x(i)=0.\tag{$15_k$}
$$
Note that the coefficients $c_{ki}$ here depend only on $\pi$,
not on $x$. Note also that $c_{kd_k}\neq0$, bacause the rational
linear relation from which we got \thetag{$15_k$} really
involved $x(d_k)$.

We define an endomorphism $C$ of \p\ by setting
$$
C(x)(k)=\sum_{i=0}^{d_k} c_{ki}x(i)
$$
for all $x\in\p$.  Regarding the elements of \p\ as infinite
column vectors, this endomorphism is given by left
multiplication by the infinite matrix $C=(c_{ni})$.
If $x\in G$, then \thetag{$15_k$} holds, and therefore
$C(x)(k)=0$, for all but finitely many $k$, i.e., $C(x)\in\s$.  So
$C$ is an endomorphism of \p\ mapping $G$ into the free group
\s.  As $G$ is weakly \s-binding, $C$ must map \s\ into a group
of finite rank.  But $C(\s)$ includes the elements $C(e_n)$, the
columns of the matrix $C$; among these are the columns
indexed by the elements $d_k$ of $D_\pi$.  The $k$th of these
columns has a non-zero entry $c_{kd_k}$ in the $k$th row and
zero entries in all earlier rows.  Thus, these columns by
themselves form a lower triangular matrix with non-zero
diagonal entries.  They are therefore linearly independent,
contrary to the fact that they lie in a group of finite rank.

\therosteritem2$\geq$\therosteritem3:\enspace trivial.

\therosteritem3$\geq$\therosteritem1:\enspace
Suppose $\Cal E$ is a family of functions $\omega\to\z$ such
that every forgetful linear predictor is evaded by some element
of $\Cal E$.  So $\Cal E$ is a subset of \p; let $G$ be the
subgroup that it and \s\ generate.  $G$ has the same cardinality
as $\Cal E$, so to complete the proof it suffices to show that $G$
is weakly \s-binding.

Suppose it were not.  Fix an endomorphism $f$ of \p\ mapping
$G$ into \s\ and mapping \s\ onto a group of infinite rank.
Each component of $f$, mapping \p\ to \z, has, by Specker's
theorem [19 Satz III], the form
$f(x)(n)=\sum_i c_{ni}x(i)$, where for each fixed $n$ only
finitely many $c_{ni}$ are non-zero.  Thus, in the matrix
$C=(c_{ni})$, each row has only finitely many non-zero entries.
So does each column, for the $i$th column is $f(e_i)\in
f(G)\subseteq\s$. (Recall that $G$ was defined so as to contain
all the $e_i$.)

We inductively choose infinitely many rows, the rows indexed
by $i_0$, $i_1$, etc., as follows.  Choose $i_0$ so that the $i_0$th
row of $C$ isn't zero.  (This is possible as $f$ is not identically
zero.)   For the induction step, suppose we have already chosen
$i_k$, an index of a non-zero row in $C$.  Let $d_k$ be the
number of the column in which the last non-zero entry of row
$i_k$ occurs.  Then choose $i_{k+1}$ so that row $i_{k+1}$ of $C$
is non-zero, but its entries in columns 0 through $d_k$ are all
zero.  The first half of this constraint is satisfied by infinitely
many rows, for otherwise $f(\p)$ would have finite rank,
whereas we are assuming that even $f(\s)$ has infinite rank.
The second half of the constraint is satisfied by all but finitely
many rows, because each column of $C$ has only finitely many
non-zero entries.  So the choice of $i_{k+1}$ is possible.

Define a forgetful linear predictor $\pi$ by $D_\pi=\{d_k\mid
k\in\omega\}$ (where $d_k$ is, as in the preceding paragraph,
the largest $d$ for which $c_{i_kd}\neq0$) and
$$
\pi_{d_k}(t)=-\frac{1}{c_{i_kd_k}}\sum_{r<d_k}c_{i_kr}t(r).
$$
This is clearly a linear predictor; it is forgetful because
$c_{i_{k+1}r}=0$ for $r\leq d_k$ by the second half of the
constraint on the choice of $i_{k+1}$.

As $G$ includes $\Cal E$ which has elements to evade any given
forgetful linear predictor, fix an $x\in G$ evading this $\pi$.
This means that, for infinitely many $k$,
$$
x(d_k)\neq-\frac{1}{c_{i_kd_k}}\sum_{r<d_k}c_{i_kr}x(r).
$$
Clearing denominators and transposing the negative terms, we
find that the inner product of $x$ and the $i_k$th row of $C$ is
non-zero.  This inner product is the $i_k$th component of $f(x)$,
so its being non-zero for infinitely many $k$ means that
$f(x)\notin\s$.  But this is a contradiction because $x\in G$ and
$f(G)\subseteq\s$.
\qed\enddemo

\proclaim{Corollary 11}
$\pp\leq\ee_l$
\endproclaim

\demo{Proof}
Theorems 9 and 10.
\qed\enddemo

By Theorem 10, any lower bound for $\ee_l=\ee_{fl}$ is also a
lower bound for all the cardinals except $\aleph_1$ listed in
Corollary~8, in particular for $\se$.  Apart from $\pp$, we have
one more such bound, given by the following theorem.

\proclaim{Theorem  12}
$\add(L)\leq\ee_l$.
\endproclaim

\demo{Proof}
We use the combinatorial description of $\add(L)$, at the end of
Section 1, as the smallest number of of functions that do not all
go through a single slalom.  We observe that, although that
description is phrased in terms of functions whose values are
elements of $\omega$ and slaloms whose values are subsets of
$\omega$, it would make no difference if $\omega$ were
replaced in both places with some other countably infinite set
$A$.  (The domains of the slaloms and of the functions under
consideration are still $\omega$; only the ranges are modified.
Thus, the requirement in the definition of slalom that
$|s(n)|=(n+1)^2$ still makes sense.)

We begin by describing a convenient $A$ for this proof.
Partition $\omega$ into finite blocks of consecutive numbers,
say $I_0=[0,a_1)$, $I_1=[a_1,a_2)$, etc., in such a way that
$I_n$ has more than $(n+1)^2$ elements.  (For example, one
could take $a_n=2n^3$.)  Let $A$ be the set of all functions with
domain equal to one of the $I_n$'s and with values in \z.  As
$A$ is countably infinite, the remarks in the preceding
paragraph apply to it.  The proof will be complete if we find
$\ee_l$ functions $\omega\to A$ that do not all go through any
single slalom $s$ (whose values are subsets of $A$).

Let $\Cal E$ be a family of $\ee_l$ functions $\omega\to\z$ as
in the definition of $\ee_l$, i.e., every linear predictor is evaded
by some $f\in\Cal E$.  For each $f\in\Cal E$, define
$f':\omega\to A$ by $f'(n)=f\rt I_n$.  Suppose $s$ were a slalom
that all these $f'$'s go through; we shall complete the proof by
deducing a contradiction.

Temporarily fix some $n\in\omega$.  Let $T$ be the set of
functions $I_n\to\z$ that are elements of $s(n)$. ($s(n)$ may
also have elements with domain $I_m$ for $m\neq n$, but
these are irrelevant to our purpose as well as to the assumption
that all $f'$ go through~$s$.)  $T$ consists of at most
$|s(n)|=(n+1)^2$ elements of the rational vector space of
rational-valued functions on $I_n$.  This vector space has
dimension $|I_n|>(n+1)^2$, so there is a non-zero linear
functional $\phi_n$ on this vector space annihilating all the
elements of $T$.  It has the form $\phi_n(t)=\sum_{i\in
I_n}c_{ni}t(i)$, for some rational coefficients $c_{ni}$, not all
zero.  Let $d_n$ be the largest $i$ for which $c_{ni}\neq0$.
We record for future use the fact that, if $f'(n)\in s(n)$ then
$\phi_n(f\rt I_n)=0$. This is because $f\rt I_n=f'(n)$ belongs to
$s(n)$ and maps $I_n$ into \z\ and therefore belongs to $T$.

Now unfix $n$, but keep the notations $\phi_n$, $c_{ni}$, and
$d_n$ from the preceding paragraph. (The subscripts $n$,
superfluous before, now prevent ambiguity.) Define a predictor
$\pi$ by setting $D_\pi=\{d_n\mid n\in\omega\}$ and, for each
$d_n\in D_\pi$, setting
$$
\pi_{d_n}(t)=-\frac{1}{c_{nd_n}}
\sum\Sb i<d_n\\i\in I_n\endSb c_{ni}t(i).
$$
This is a forgetful linear predictor.  If some $f:\omega\to\z$ and
some $n\in\omega$ satisfy $\phi_n(f\rt I_n)=0$, then trivial
algebraic manipulation of this equation gives
$f(d_n)=\pi_{d_n}(f\rt d_n)$.

For each $f\in\Cal E$, the assumption that $f'$ goes through $s$
means that for all but finitely many $n$ we have $f'(n)\in s(n)$,
and so $\phi_n(f\rt I_n)=0$, and so $f(d_n)=\pi_{d_n}(f\rt
d_n)$.  But that means that $\pi$ predicts every $f\in\Cal E$,
contrary to the choice of $\Cal E$.
\qed\enddemo

Upper bounds for the evasion cardinals $\ee$ and $\ee_l$ are less
interesting than lower bounds, since they do not imply corresponding
bounds for the more natural group-theoretic cardinals in Corollary~8.
Nevertheless, we record for the sake of completeness the upper bounds
which we have been able to obtain.

\proclaim{Theorem  13}
All three of $\unif(L)$, $\unif(B)$, and $\dd$ are upper bounds
for $\ee$.  Furthermore, $\min\{\ee,\bb\}\leq\add(B)$.
\endproclaim

\demo{Proof}
$\ee\leq\unif$:\enspace
We handle the measure and category cases together.
It is easy to verify that, for any predictor $\pi$, the set of
functions $x\in2^\omega$ (i.e., functions on $\omega$ with only
0 and 1 as values) that are predicted by $\pi$ is a set of first
category and measure zero in $2^\omega$.  The desired
inequalities follow immediately.

$\ee\leq\dd$:\enspace
Let $\Cal D$ be a family of $\dd$ functions from
$\omega\times\omega$ (the set of pairs of natural numbers) to
$\omega$ such that every function
$\omega\times\omega\to\omega$ is majorized (everywhere)
by a function from $\Cal D$. (See the remarks in Section~1 after
the definition of $\dd$.)  To each function $g\in\Cal D$, we
associate a function $x:\omega\to\omega$ by the recursion
$$
x(n)=g(n,1+\max\{x(p)\mid p<n\}).\tag{16}
$$
We shall show that every $\omega$-valued predictor $\pi$ is
evaded by one of these $\dd$ functions $x$.

Let $\pi=(D_\pi,(\pi_n))$ be given, and define a function
$f:\omega\times\omega\to\omega$ by
$$
f(n,k)=\cases
\max\{\pi_n(t)\mid t\in\omega^n\text{ and all values of
}t\text{ are }<k\},&\text{if }n\in D_\pi\\
\text{arbitrary},&\text{otherwise.}
\endcases
$$
Choose $g\in\Cal D$ such that $g(n,k)>f(n,k)$ for all $n$ and $k$,
and let $x$ be the function obtained from $g$ by \thetag{16}.
We shall show that this $x$ evades $\pi$.

Consider any $n\in D_\pi$.  Let $k=1+\max\{x(p)\mid p<n\}$
and note that $x\rt n$, being in $\omega^n$ and having all its
values $<k$, is one of the $t$'s involved in the definition of
$f(n,k)$.  So $f(n,k)\geq\pi_n(x\rt n)$.  On the other hand, by
the definition \thetag{16} of $x$ and the choice of $g$, we also
have $x(n)\geq g(n,k)>f(n,k)$.  Comparing these inequalities, we
find that $x(n)\neq\pi_n(x\rt n)$. As $n$ was an arbitrary
element of $D_\pi$, this shows that $x$ evades $\pi$. (It
actually shows more: not only are infinitely many of $\pi$'s
predictions wrong, as required for evasion, all of them are
wrong.)

$\min\{\ee,\bb\}\leq\add(B)$:\enspace
Since $\add(B)=\min\{\cov(B),\bb\}$, it suffices to prove that
$\min\{\ee,\bb\}\leq\cov(B)$.  For this we use the
combinatorial description of $\cov(B)$ at the end of Section~1.
We must show that, if $\Cal F$ is a family of functions
$\omega\to\omega$ of cardinality smaller than
$\min\{\ee,\bb\}$, then there exists $g:\omega\to\omega$ such
that every $x\in\Cal F$ is infinitely often equal to $g$.

Let such an $\Cal F$ be given.  As its cardinality is smaller than
$\bb$, fix $d:\omega\to\omega$ eventually majorizing (i.e.,
$\geq^*$) every $x\in\Cal F$.  For any natural number $k$ we
will write $\max(d,k)$ for the function whose value at any $n$
is $\max\{d(n),k\}$.  For any $x\in\Cal F$, we can find $k$ so
large that $x$ is majorized everywhere by $\max(d,k)$. ($d$
majorizes all but finitely many values of $x$, so we just choose
$k$ to majorize those finitely many values.)

Partition $\omega$ into finite blocks of consecutive numbers,
say $I_0=[0,a_1)$, $I_1=[a_1,a_2)$, etc., in such a way that the
cardinality of $I_k$ is
$$
|I_k|=\prod_{j<a_k}\max\{d(j),k\},\tag{17}
$$
i.e., the number of functions $I_0\cup\dots\cup I_{k-1} \to
\omega$ that are majorized on their domain by $\max(d,k)$.
Note that none of these blocks are empty.

Let $A$ be the set of functions into $\omega$ whose domains
are blocks in this partition.  For each $x\in\Cal F$, define an
associated function $x':\omega\to A$ by $x'(k)=x\rt I_k$.
As $|\Cal F|<\ee$, all these functions $x'$ are predicted by a
single $A$-valued predictor $\pi=(D_\pi, (\pi_n))$.

We use $\pi$ to define the desired $g$ (infinitely often equal to
each $x\in\Cal F$) as follows.  If $k\notin D_\pi$, define $g\rt
I_k$ arbitrarily.  If $k\in D_\pi$, then proceed as follows.  List
the functions $I_0\cup\dots\cup I_{k-1}\to\omega$ that are
majorized (on their domain) by $\max(d,k)$ as $t_0$, $t_1$,
\dots $t_{r-1}$; the number of such functions, here called $r$, is
the cardinality of $I_k$, by \thetag{17}.  For each $j<r$, obtain
$t'_j$ from $t_j$ analogously to the way $x'$ was defined from
$x$, i.e., $t'_j:\{0,1,\dots,k-1\}\to A$ and $t'_j(i)=t_j\rt I_i$.
Consider the set $X$ of those $\pi_k(t'_j)$ (for our fixed $k$ and
arbitrary $j<k$) that are functions $I_k\to\omega$. This is a set
of at most $r$ functions on a set $I_k$ of size $r$.  So there is a
function $I_k\to\omega$ that agrees at least once with each of
the functions in $X$.  Let $g\rt I_k$ be such a function.

To show that $g$ has the desired property, consider any
$x\in\Cal F$; we must show that $g$ and $x$ are infinitely often
equal.  Temporarily fix some $k\in D_\pi$ so large that
$\max(d,k)$  majorizes $x$ everywhere and so large that $\pi$
correctly predicts $x'(k)$, i.e., $x'(k)=\pi_k(x'\rt k)$.   Then, in
the definition of $g\rt I_k$, one of the $t_j$'s under
consideration was $x\rt I_0\cup\dots\cup I_{k-1}$ and the
corresponding $t'_j$ was $x'\rt k$.  So, as $\pi$ correctly
predicted $x'(k)$, we have that
$$
x\rt I_k=x'(k)=\pi_k(x'\rt k)=\pi_k(t'_j)\in X
$$
is one of the functions with which we made $g\rt I_k$ agree at
least once.  This proves that $g$ and $x$ agree at least once in
$I_k$.  Now un-fix $k$.  The preceding argument applies to
infinitely many $k$'s, namely all sufficiently large elements of
$D_\pi$.  So $g$ and $x$ are infinitely often equal.
\qed\enddemo

\proclaim{Corollary  14}
$\ee_l\leq\add(B)$
\endproclaim

\demo{Proof}
Combine Theorem~13 with the facts that $\ee_l\leq\ee$ (as we
remarked right after defining these cardinals) and
$\ee_l\leq\bb$, which follows from Theorem~10 and
Corollary~8.
\qed\enddemo

\remark{Remark}
Referring to the intuitive interpretation of predictors, which
was described just after their definition, we point out that there
is a natural extension of the concept of predictor that matches
even better the intuition of a guessing strategy.  Instead of
having a fixed set $D_\pi$ of numbers $n$ for which the
strategy tries to predict the value $x(n)$ of an unknown
function, given the list $x\rt n$ of prior values, we could let the
decision whether to attempt a prediction of $x(n)$ for a
particular $n$ depend also on $x\rt n$. Of course, we should
either require that, for each $x$, the strategy attempts infinitely
many predictions or else declare that any $x$ for which this
fails is deemed to have evaded $\pi$.  The cardinal $\ee+$
associated to this concept, namely the smallest number of
functions $\omega\to\omega$ needed to evade all
$\omega$-valued predictors of this generalized sort, is clearly
$\geq\ee$.  All the upper bounds we have given for $\ee$
apply also to $\ee+$, with only minor modifications of the
proofs.

We also remark that, if we were to weaken the definition of
``predicts'' by requiring only that infinitely many (rather than
all but finitely many) predictions are correct, then the new
$\ee$ would be $\geq$ the old one, but it would still be
$\leq\dd$ and its minimum with $\bb$ would still be
$\leq\add(B)$, by the same proofs as for the original $\ee$.
\endremark

%% file: ccpref.tex
%
